\documentclass{article}
\usepackage{hyperref}
\usepackage{graphicx}
\usepackage{amssymb,amsmath}
\usepackage{color}
\newcommand{\footremember}[2]{%
    \footnote{#2}
    \newcounter{#1}
    \setcounter{#1}{\value{footnote}}%
}
 
\author{Franti\v sek Kardo\v s\footremember{labri}{Univ. Bordeaux, CNRS, Bordeaux INP, LaBRI, UMR 5800, Talence, France}%
\and Martina Mockov\v ciakov\'a\footremember{plzen}{University of West Bohemia, Pilsen, Czech Republic}%
}
\title{At least half of the leapfrog fullerene graphs have exponentially many Hamilton cycles}

\newtheorem{theorem}{Theorem}
\newtheorem{lemma}{Lemma}
\newtheorem{claim}{Claim}

\newcommand{\ep}{\hfill $\square$}

\begin{document}

\maketitle

\begin{abstract}
A fullerene graph is a 3-connected cubic planar graph with pentagonal and hexagonal faces. The leapfrog transformation of a planar graph produces the trucation of the dual of the given graph. A fullerene graph is leapfrog if it can be obtained from another fullerene graph by the leapfrog transformation. We prove that leapfrog fullerene graphs on $n=12k-6$ vertices have at least $2^{k}$ Hamilton cycles.
\end{abstract}

\section{Introduction}

Tait conjectured that every 3-connected cubic planar graph is hamiltonian. Had his conjecture been true, it would have implied the Four Colour Theorem.
This conjecture was disproved by Tutte in 1946. Tutte also proved (1956) that every 4-connected planar graph is hamiltonian \cite{tutte}. 

However, every known non-hamiltonian 3-connected cubic planar graph contains a face of size at least seven. Barnette conjectured in 1969 (and Goodey stated it in an informal way as well), that all 3-connected cubic planar graphs with faces of size at most 6 are hamiltonian. 

In 1975 Goodey proved that all 3-connected cubic planar graphs with faces of size 4 and 6 are hamiltonian\cite{good1}, and in 1977 he proved that it holds also for such graphs with faces of size 3 and 6 \cite{good2}.
 All 3-connected cubic planar graphs with faces of size at most 6 have recently been proved to be hamiltonian by Kardos \cite{ferdo}, confirming Barnette-Goodey conjecture.

A \emph{fullerene graph} is a 3-connected cubic planar graph with pentagonal and hexagonal faces.
Fullerene graphs are of a standalone interest since they are used to modelize all-carbon molecules of spherical shape -- carbon atoms are represented by the vertices and bonds between adjacent atoms are represented by the edges of the graph. 

A fullerene is called a \emph{leapfrog fullerene}, if it can be constructed from other fullerene graph $G$ by so-called leapfrog transformation -- the truncation of the dual of $G$. We will use the notation $H=L(G)$ whenever $H$ is a leapfrog of $G$. If this is the case, then there are three types of faces in $H$: pentagonal faces of $H$ correspond to pentagonal faces of $G$, some hexagonal faces correspond to hexagonal faces of $G$, and the remaining hexagonal faces of $H$ correspond to vertices of $G$. The set of the faces of the first two types (those that correspond to the faces of $G$) is independent: their facial cycles form a 2-factor of $H$.
It is straightforward to see that $|V(L(G))|=3|V(G)|$.

\medskip
Maru\v si\v c \cite{mar} proved that every leapfrog fullerene graph obtained from a fullerene graph with an odd number of faces is hamiltonian.

Both Maru\v si\v c's and Kardo\v s's proofs rely on the correspondence between Hamilton cycles in a cubic planar graph and stable-tree decompositions in a suitably defined residual graph.

\begin{figure}
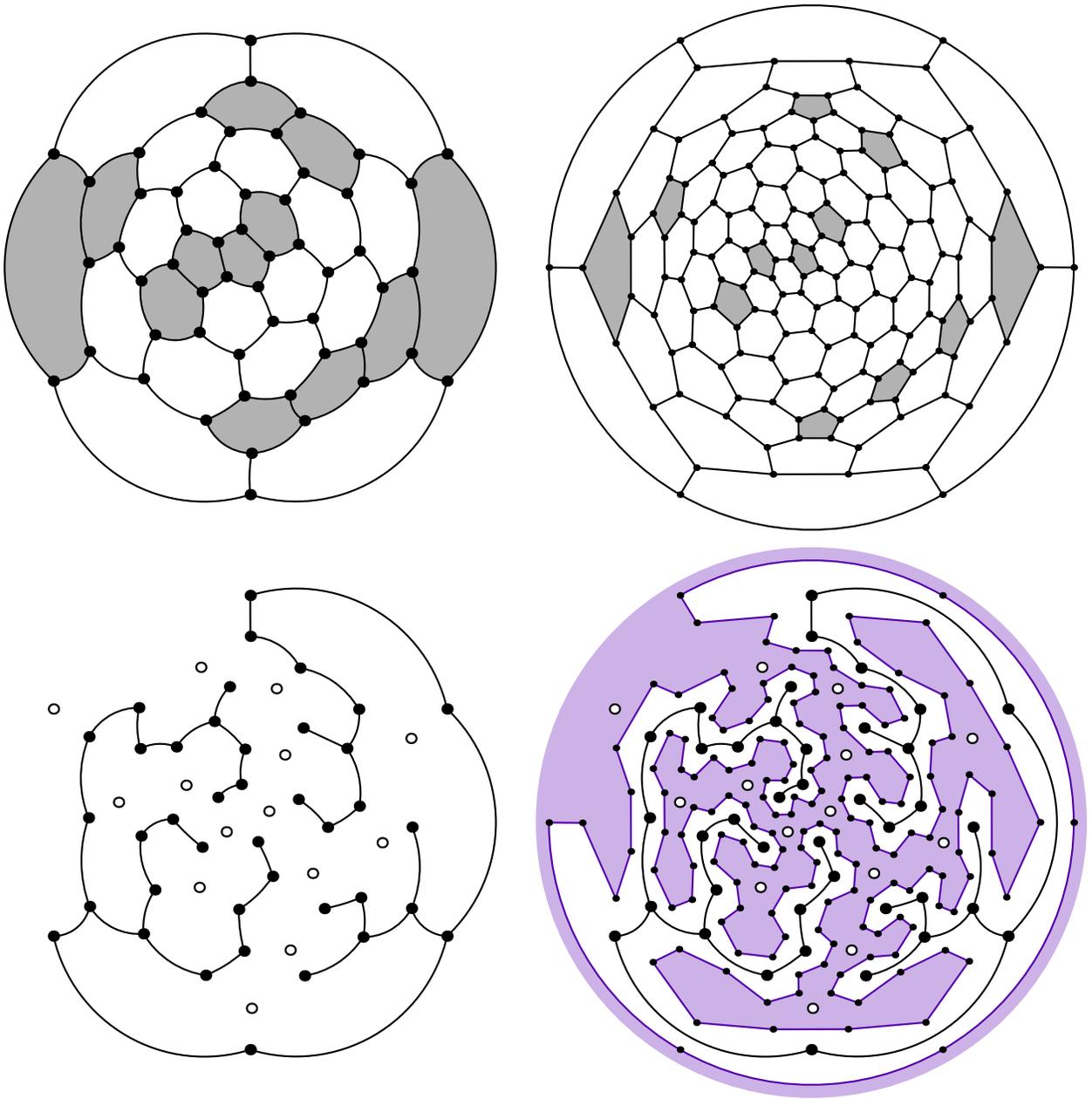

\centerline{
\includegraphics[scale=1]{mar.1}
\hfil
\includegraphics[scale=1]{mar.0}
}
\centerline{
\includegraphics[scale=1]{mar.41}
\hfil
\includegraphics[scale=1]{mar.42}
}
\caption{An example of a fullerene graph $G$ and the fullerene graph $H$ obtained from $G$ by leapfrog transformation (top, left to right). A stable-tree decomposition of $G$ and the corresponding Hamilton cycle in $H$ (bottom, left to right).}
\label{fig:ex}
\end{figure}

A \emph{stable-tree decomposition} of a graph is a partition of its vertex set into two sets, say $W$ and $B$, such that $W$ is a stable set and $B$ induces a tree. We call the vertices in $W$ white and those in $B$ black. If $H=L(G)$, then every stable-tree decomposition of $G$ corresponds to a Hamilton cycle in $H$, see Figure \ref{fig:ex} for illustration.

To prove Maru\v si\v c's result, it suffices therefore to observe that for fullerene graphs with an odd number of faces the existence of a stable-tree decomposition is guaranteed by

\begin{theorem}[Payan and Sakarovitch, 1975 \cite{PS}]
Let $G$ be a cubic graph on $n=4k-2$ vertices ($k\ge 1$). If $G$ is cyclically $4$-edge-connected, then $H$ has a stable-tree decomposition.
\label{th:pyber}
\end{theorem}
and
\begin{theorem}[Do\v sli\' c, 2003 \cite{doslic}]
Let $G$ be a fullerene graph. Then $G$ is cyclically 5-edge-connected.
\label{th:doslic}
\end{theorem}

In this paper we prove

\begin{theorem}
Let $H$ be a leapfrog fullerene graph of a fullerene graph $G$ on $n=4k-2$ vertices. Then $H$ has at least $2^{k}$ Hamilton cycles.
\label{th:main}
\end{theorem}

\section{Preliminaries}

In this section we introduce some technical notions used in the proof of Theorem \ref{th:main}.

An \emph{improper stable-tree decomposition} of a cubic planar graph $G$ is a partition of the vertex set of $G$ into two sets, say $W$ and $B$, such that $W$ is a stable set, $B$ induces a forest with exactly three components, and there exists a hexagonal face of $G$ incident to one vertex of each component of $G[B]$.

In this context, a hexagonal face of $G$ incident to one vertex of each of the three components of $G[B]$ will be called a \emph{graceful} hexagon. A vertex in $W$ adjacent to one vertex of each of the three components of $G[B]$ will be called a \emph{graceful} vertex.

A \emph{generalized stable-tree decomposition} of a cubic planar graph $G$ is either a (proper) stable-tree decomposition of $G$ or an improper stable-tree decomposition of $G$.

\begin{lemma}
Let $(W,B)$ be an improper stable-tree decomposition of a fullerene graph $G$, let $f$ be a graceful hexagon in $G$. Then either there exists another graceful hexagon in $G$ distinct from $f$, or there exists a graceful vertex in $G$.
\end{lemma}

Proof. Let $B_1$, $B_2$, $B_3$ be the three components of $G[B]$. 
Let $F_1$ be the set of faces of $G$ incident to a vertex of $B_1$ and at least one vertex of at least one out of $B_2$ and $B_3$. 
Clearly $f\in F_1$. 

If there is another face $f'$ in $F_1$ incident to a vertex both in $B_2$ and $B_3$, the three vertices from the three different components have to be separated by white vertices on the facial cycle of $f'$. Since in $G$ the maximum face size is 6, there is exactly one vertex from each component, and so $f'$ is a graceful hexagon.

Suppose now that there are no graceful hexagons other that $f$ in $G$. 
Tracing the faces in $F_1$ in a cyclic order around the component $B_1$ one can find two adjacent faces, say $f_2$ and $f_3$, incident each to a vertex of $B_2$ and $B_3$, respectivement. Then the edge they share has at least one black vertex, say $v_1$, which has to be from $B_1$, and a white vertex, say $v$, which has three black neighbors, say $v_1$, $v_2$, and $v_3$, with $v_i$ incident to $f_i$ for $i=2,3$. It is easy to see that $v_i \in B_i$ for $i=2,3$, otherwise there could not be any graceful hexagons in $G$. \ep

\begin{lemma}
Let $(W,B)$ be a proper stable-tree decomposition of a cubic graph $G$ on $n$ vertices. Then $n=4k-2$ for some $k\in\mathbb{N}$, moreover, $|W|=k$. 
\end{lemma}

Proof. Let $b=|B|$ and $w=|W|$. Since the black vertices induce a tree, the number of edges joining two black vertices is $b-1$. Since the white vertices are stable, the number of edges joining a black and a white vertex is $3w$. Since the graph is cubic, the overall number of edges is $3n/2$. Therefore,
$$
\gathered
b-1 + 3w = \frac{3(b+w)}{2} \\
3w - 2 = b
\endgathered
$$
and so $n=b+w=4w-2$; it suffices to set $k=w$.
\ep

\begin{lemma}
Let $(W,B)$ be an improper stable-tree decomposition of a cubic graph $G$ on $n$ vertices. Then $n=4k-2$ for some $k\in\mathbb{N}$, moreover, $|W|=k+1$. 
\end{lemma}

Proof. Let $b=|B|$ and $w=|W|$. Since the black vertices induce three distinct trees, the number of edges joining two black vertices is $b-3$. Since the white vertices are stable, the number of edges joining a black and a white vertex is $3w$. Since the graph is cubic, the overall number of edges is $3n/2$. Therefore,
$$
\gathered
b-3 + 3w = \frac{3(b+w)}{2} \\
3w - 6 = b
\endgathered
$$
and so $n=b+w=4w-6 = 4 (w-1)-2$; it suffices to set $k=w-1$.
\ep

\bigskip
Let $G$ be a 2-connected graph. An \emph{ear decomposition} of $G$ is a partition $(E_0,E_1,\dots,E_t)$ of the set of edges of $G$ such that:
\begin{itemize}
\item $E_0$ is a cycle;
\item $E_k$ is a path edge-disjoint from $G_{k-1}:=\cup_{i=0}^{k-1} E_i$ with two distinct end-vertices both in $G_{k-1}$.
\end{itemize}

A vertex contained in $E_i$ but not in $G_{k-1}$ is an \emph{internal} vertex of the ear $E_i$.

It is known (add citation Whitney) that every 2-connected graph admits an ear-decomposition.

Let $G$ be a 3-connected planar graph. Denote $G^*$ the dual graph of $G$.
A \emph{dual search order} is an ordering $(f_0,f_1,\dots,f_s)$ of faces of $G$ such that $f_k$ is adjacent to at least one from $f_0,\dots,f_{k-1}$. Typically, a BFS, DFS or LexBFS of $G^*$ provides a dual search order.

A \emph{facial ear decomposition} of $G$ is an ear decomposition of $G$ corresponding to a dual search order: $E_0$ is the facial cycle of $f_0$; for each face $f_k$, the difference between $G[f_0,\dots,f_{k-1}]$ and $G[f_0,\dots,f_{k-1},f_k]$ is an edge-disjoint union (possibly empty) of ears with the endvertices in the former.

An ear decomposition of a 3-connected planar graph is \emph{nice} if it is a facial ear decomposition of $G$ such that
\begin{itemize}
\item $f_0$ is a hexagon (and hence $E_0$ is a 6-cycle)
\item each $E_j$ ($j\ge 1$) is a path on at most 4 edges.
\end{itemize}

\begin{lemma}
Let $G$ be a fullerene graph. Then there exists a nice ear decomposition of $G$.
\end{lemma}

Proof. It suffices to consider a LexBFS of $G^*$ with the first vertex corresponding to a hexagon and the second to a pentagon. Then each new vertex is adjacent to at least two vertices already visited, so the corresponding ear(s) is (are) composed of at most four edges. \ep

\section{Results}

We are now ready to prove the main theorem. In order to proof the existence of exponentially many Hamilton cycles in a leapfrog fullerene graph, we find exponentially many generalized stable-tree decompositions in the original/underlying fullerene graph, using the following invariant: Every time a white vertex (a vertex in the stable set) is to be introduced, we can do it in two different ways; in total, there is a linear number of white vertices in any generalized stable-tree decomposition.

Proof of Theorem \ref{th:main}. We will introduce a procedure that produces a set of  generalized stable-tree decompositions of $G$. Every proper stable-tree decomposition of $G$ will then correspond to a distinct Hamilton cycle of $H$, every improper stable-tree decomposition of $G$ will correspond to two distinct Hamilton cycles of $H$.

Let $(E_0,E_1,\dots,E_t)$ be a nice ear decomposition of $G$.

Let $G^{j\to} := G[\cup_{i=j}^t V(E_i)]$ be the graph induced on the set of internal vertices of ears starting form the $j$-th ear. Clearly $G^{j\to}$ is an empty graph for $j>t$ and $G^{0\to}=G$, since every vertex of $G$ is an internal vertex of exactly one ear.

In the next paragraphs we define inductively how to decide for vertices of $G$ whether they are black or white. We decide about the color of internal vertices of the ears in the descending order, in such a way that after each step, the set of the colored vertices is equal to $G^{j\to}$ for some $j$.

Let $(W,B)$ be a decomposition of $G^{j\to}$. A black vertex $v$ of $G^{j\to}$ is called a \emph{contact vertex} if it is adjacent to a vertex in $G\setminus G^{j\to}$.

We prove that for each $j$ there exists a set $\mathcal{D}_j$ of decompositions $(W^j_i,B^j_i)$ of $G^{j\to}$ with the following properties:
\begin{itemize}
\item[(i)] $W^j_i$ is a stable set in $G^{j\to}$ as well as in $G$,
\item[(ii)] $B^j_i$ induces a forest in $G^{j\to}$, moreover, every component of $G^{j\to}[B^j_i]$ contains a contact vertex.
\end{itemize}

In other words, we keep the property that every black component is a tree attached to at least one vertex of some non-colored ear.

In order to enumerate $\mathcal{D}_1$, we introduce a parent-child relation between the decompositions in $\mathcal{D}_{j+1}$ and $\mathcal{D}_{j}$, which is defined together with the decompositions themselves by induction in the following way:

For $j>t$, we define $\mathcal{D}_j:=\{(\emptyset,\emptyset)\}$, we call this decomposition \emph{trivial}.

(Initialisation) Let $j_{max}$ be the largest $j$ such that $E_j$ has at least two edges (at least one internal verex). We color all the internal vertices of $E_j$ black, i.e. we set $\mathcal{D}_{j_{max}}=\{(\emptyset,int(E_{j_{max}}))\}$, where $int(E_{j_{max}})$ is the set of internal vertices of $E_{j_{max}}$. Clearly, the black vertices induce a forest (a path), while the set of white vertices is stable. We call this decomposition \emph{root decomposition}, it is the only non-trivial decomposition without a parent.

(Propagation) Let $E_j$ be an ear. We distinguish cases according to the number of edges of $E_j$ (number of internal vertices of $E_j$). 

CASE 0: If the ear $E_j$ only has one edge, then it has no internal vertices, thus $G^{j\to}=G^{j+1\to}$ and so we set $\mathcal{D}_{j}:=\mathcal{D}_{j+1}$; we do not distinguish between the two sets of decompositions.

CASE 1: Let $E_j$ be an ear with two edges and one internal vertex $v$. Then $v$ has exactly one neighbor in $G^{j+1\to}$, say $u$. We set
$$
\mathcal{D}_j=\{(W,B\cup \{v\}) | (W,B)\in \mathcal{D}_{j+1}\}.
$$
Informally, for every decomposition from $\mathcal{D}_{j+1}$, we color the vertex $v$ black. Let us check that $\mathcal{D}_j$ satisfies (i) and (ii).
Clearly, $W$ remains a stable set. If $u$ is white then $v$ becomes a new component in the forest induced on the black vertices in $G^{j\to}$; if $u$ is black then $v$ is added to an existing component as a new leaf. In both cases $v$ is a contact vertex of the component it is contained in; all other black components keep their contact vertices.

We say that the decomposition $(W,B\cup\{v\})\in \mathcal{D}_j$ of $G^{j\to}$ is (the only) child of $(W,B)\in \mathcal{D}_{j+1}$. It is easy to see that the parent-child relation between $\mathcal{D}_{j+1}$ and $\mathcal{D}_j$ is a bijection.

CASE 2: Let $E_j$ be an ear with three edges and two internal vertices $v_1$ and $v_2$. Each of the vertices $v_1$ and $v_2$ has exactly one neighbor in $G^{j+1\to}$, say $u_1$ and $u_2$, respectively. For every decomposition in $\mathcal{D}_{j+1}$ the vertices $u_1$ and $u_2$ are already colored.

Let $(W,B)$ be a decomposition of $G^{j+1\to}$. We consider subcases with respect to the color of $u_1$ and $u_2$.

If both $u_1$ and $u_2$ are white, then we color $v_1$ and $v_2$ both black, introducing a new black component having two vertices, both being contact vertices.

If one of $u_1$ and $u_2$, say $u_1$, is white and the other one, say $u_2$, is black, then we color $v_1$ and $v_2$ both black, adding a path on two edges to an existing black component of $B$. Both $v_1$ and $v_2$ are contact vertices of this component.

It remains to consider the case when both $u_1$ and $u_2$ are black. If they belong to two different black components, then, similarly to the previous cases, we color $v_1$ and $v_2$ black, merging two black components into one with $v_1$ and $v_2$ as contact vertices. Up to this point we always color the two new vertices black, and so we define a unique child in $\mathcal{D}_j$ for a given decomposition from $\mathcal{D}_{j+1}$.

The last case (and the only interesting one) is when $u_1$ and $u_2$ belong to the same black component. This is the first case when for a decomposition from $\mathcal{D}_{j+1}$ we define two children in $\mathcal{D}_j$: We color $v_1$ and $v_2$ black and white (one vertex per color), in both possible ways. Since one of them is white, we do not create a cycle in the component containing them, and since the other one is black, the component still contains at least one contact vertex. 

CASE 3: Let $E_j$ be an ear with four edges and three internal vertices $v_1$, $v_2$, and $v_3$. Each of them has exactly one neighbor in $G^{j+1\to}$, say 
$u_1$, $u_2$, and $u_3$,
respectively. 

If each $u_i$ is either white or belongs to a distinct black component, we color all of $v_1$, $v_2$, $v_3$ black, eventually merging some existing black components into a new black component having $v_1$ and $v_3$ as contact vertices. No new white vertex is created, and a unique child is defined.

\begin{figure}[ht]
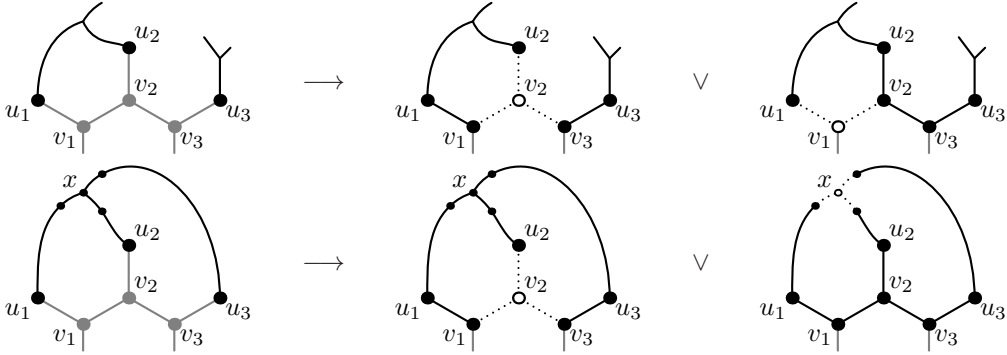

\centerline{
\begin{tabular}{ccccc}
\begin{tabular}{c}
\includegraphics[scale=1]{cases.30}
\end{tabular}
&$\longrightarrow$&
\begin{tabular}{c}
\includegraphics[scale=1]{cases.31}
\end{tabular}
&$ \vee $&
\begin{tabular}{c}
\includegraphics[scale=1]{cases.32}
\end{tabular}
\\
\begin{tabular}{c}
\includegraphics[scale=1]{cases.33}
\end{tabular}
&$\longrightarrow$&
\begin{tabular}{c}
\includegraphics[scale=1]{cases.34}
\end{tabular}
&$ \vee $&
\begin{tabular}{c}
\includegraphics[scale=1]{cases.35}
\end{tabular}
\end{tabular}
}
\caption{If an ear with three internal vertices is incident to two different components of the black forest, there are two choices of the position of a new white vertex (top line). If an ear with three internal vertices is incident to a single black component thrice, a second solution can be found by recoloring a black vertex to white.}
\label{fig:cases}
\end{figure}

If exactly two among $u_i$, $i=1,2,3$, belong to the same black component and the third one is either white or contained in another black component, then we color the former two with black and white (two possibilities) and the third one black. It is easy to check that we obtain a decomposition satisfying both conditions (i) and (ii). There is always one new white vertex, and two children are defined, see the first row of Figure \ref{fig:cases} for illustration.

The last remaining case is when all the three $u_i$s belong to the same black component, say $C$. Since $C$ is a tree, for each pair $u_i,u_{i'}$, $1\le i<i'\le 3$ there is a unique path from $u_i$ to $u_{i'}$ in $C$; the three paths have to have one (and only one) vertex in common, say $x$.
We introduce two children decompositions: In the first, we color $v_1$ and $v_3$ black, and $v_2$ white. In the second, we color all $v_1$, $v_2$, $v_3$ black and we recolor $x$ from black to white. 
Observe that in both cases, $u_1$, $u_2$, and $u_3$ belong to the same black component, having $v_1$ and $v_3$ as contact vertices.

\bigskip

So far, for each $j\ge 1$ and for each decomposition in $\mathcal{D}_{j+1}$ we have defined either a single child or two children in $\mathcal{D}_j$, by coloring up to three new vertices, and eventually recoloring one vertex already colored. We need to check that it is not possible to obtain the same decomposition in $\mathcal{D}_j$ from two distinct decompositions in $\mathcal{D}_{j+1}$ in two different ways.

We introduce a series of claims, each of them being either a direct consequence of the definititions of the decompositions above, or a corollary of the previous claims.

\begin{claim} 
For every $j\ge 2$ and for every decomposition $D$ in $\mathcal{D}_j$, $D$ has either one child, with the same number of white vertices, or it has two children, both with one more white vertex.
\end{claim}
\begin{claim} 
Whenever two black vertices start to belong to the same black component, they continue to belong to the same black component unless (at least) one of them is recolored white.
\end{claim}
\begin{claim} 
At the moment when an uncolored vertex becomes white, it has at least two black neighbors belonging to the same black component. Moreover, at any time, a white vertex has at least two black neighbors belonging to the same black component. 
\label{claim:3}
\end{claim}
\begin{claim} 
Every non-root decomposition has exactly one parent.
\end{claim}

Indeed, for all the cases but the last subcase of CASE 3 the parent decomposition can be restored by simply uncoloring the vertices of the last ear. If a decomposition created in the last subcase of CASE 3 by recoloring $x$ white and coloring $v_1$, $v_2$, and $v_3$ black had another parent, then it only could be the decomposition obtained by uncoloring $v_1$, $v_2$, and $v_3$ without changing the color of $x$; in that decomposition the three black neighbors of $x$ would belong to three different black components, a contradiction with Claim \ref{claim:3}.

\medskip
To complete the proof of the main theorem, it suffices to define how each of the decompositions in $\mathcal{D}_1$ extends to $E_0$. Recall that $E_0$ is a 6-cycle. Let $V(E_0)=\{v_1,\dots,v_6\}$ in a cyclic order, let $u_i$ be the neighbor of $v_i$ outside $E_0$, $i=1,\dots,6$.

Let $(W_1,B_1)$ be a decomposition of $G^{1\to}$. Let $w_0$ be the number of white vertices adjacent to vertices of $E_0$. Let $c$ be the number of components of the forest $B_1$. Let $w_1=|W_1|$, let $b_1=|B_1|$. 

\begin{claim} $w_0+c\in\{1,3,5\}$. Moreover, $w_1=k-\frac{7-(w_0+c)}{2}$.
\end{claim}

Proof. We have $w_0+c\le 6$ since each black component of $W_1$ has a black contact vertex.

The number of edges joining two black vertices in $B_1$ is equal to $b_1-c$. The number of edges joining a black vertex in $B_1$ and a white vertex in $W_1$ is equal to $3w_1-w_0$. 
On the other hand, the total number of vertices of $G^{1\to}$, $b_1+w_1 = n-6 = 4k-8$, and 
the overall number of edges of $G^{1\to}$ is 
$$b_1-c+3w_1-w_0 = |E(G)|-12 = \frac32 |V(G)|-12 = \frac32(4k-2)-12 = 6k-15.$$
By eliminating $b_1$ we get $2w_1-w_0-c=2k-7$, so $c+w_0$ is odd and the claim is now verified. \ep

\bigskip

If $w_0+c = 5$, then there is exactly one pair of black vertices $u_{i_1}$, $u_{i_2}$ belonging to the same component; every other $u_i$ is either a white vertex or belongs to a distinct component.
We color $v_{i_1}$ and $v_{i_2}$ black and white, in both possible ways; we color all other non-colored vertices black. It is easy to check that two different proper stable-tree decompositions of $G$ are obtained this way.

If $w_0+c=1$, then $w_0=0$ and $c=1$ since every black component has at least one black contact vertex. It means that all the vertices $u_1,\dots,u_6$ are black and they belong to the same black component. We define eight different proper stable-tree decompositions of $G$ in the following way:

First, we color $v_1,v_3,v_5$ black and $v_2,v_4,v_6$ white. This gives the first decomposition. Then, we choose one of the three white vertices $v_i$ and recolor it black (there are three possible choices), and search in the black tree $B_1$ for the vertex $x$ which is the intersection of the paths joining $u_{i-1}$, $u_i$, and $u_{i+1}$ (indices modulo $6$); we recolor $x$ white. This gives three other decompositions.

It follows from Claim \ref{claim:3} that the four decompositions are distinct from each other.

Next, we find four other decompositions in the same way, starting with $v_1,v_3,v_5$ white and $v_2,v_4,v_6$ black.

\begin{figure}[ht]
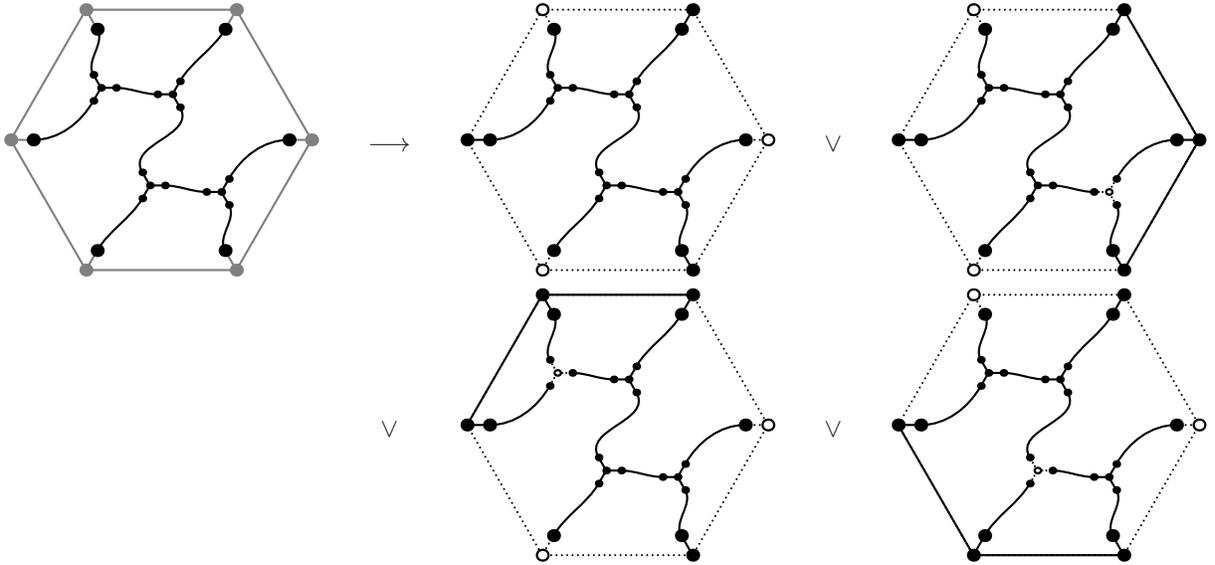

\centerline{
\begin{tabular}{ccccc}
\begin{tabular}{c}
\includegraphics[scale=1]{six.0}
\end{tabular}
&$\longrightarrow$&
\begin{tabular}{c}
\includegraphics[scale=1]{six.1}
\end{tabular}
&$ \vee $&
\begin{tabular}{c}
\includegraphics[scale=1]{six.2}
\end{tabular}
\\
&$ \vee $&
\begin{tabular}{c}
\includegraphics[scale=1]{six.3}
\end{tabular}
&$ \vee $&
\begin{tabular}{c}
\includegraphics[scale=1]{six.4}
\end{tabular}
\end{tabular}
}
\caption{If all the six vertices of the hexagon $E_0$ (depicted here as the outer face) are adjacent to black vertices belonging to the same component, we can first choose one of the two options how to them alternatively black and white, and then for every white vertex we can exchange it for a black internal vertex.}
\label{fig:six}
\end{figure}

The last case to consider is $w_0+c=3$. If we consider the numbers of $u_i$s contained in black components (a white $u_i$ being treated as a distinct (empty) component), there are three subcases: $(4,1,1)$, $(3,2,1)$, $(2,2,2)$.

Let there be a component containing four of $u_i$s. We can proceed in an analogous way as in the previous case in order to find six distinct decompositions, whereas we only need four.

Let the black components contain $3$, $2$, and $1$ of the $u_i$s. We can treat the one having 3 in the same way as an ear with three internal vertices adjacent to the same component, and the one having 2 in the same way as an ear with two internal vertices adjacent to the same component. Combining these partial resolutions gives four distinct decompositions, as needed.

The last case to consider is $(2,2,2)$. It implies $w_0=0$, there are exactly three black components, each with two contact vertices. There are two subcases. Up to symmetry, thanks to planarity, in the same component are either $u_1$ with $u_2$, $u_3$ with $u_4$, and $u_5$ with $u_6$, or $u_1$ with $u_4$, $u_2$ with $u_3$, and $u_5$ with $u_6$.

In both subcases, we color $v_1,v_3,v_5$ black and $v_2,v_4,v_6$ white or vice versa -- we introduce three white vertices and get two generalized stable-tree decompositions, each of them at the end will give two Hamilton cycles in the leapfrog fullerene graph $H$.

To conclude, we have defined a rooted binary tree of depth $k$, in which leaves at depth $k$ correspond to proper stable-tree decompositions of $G$ and leaves at depth $k-1$ correspond to improper stable-tree decompositions of $G$.
\ep

\section{Concluding remarks}

The method introduced in the paper can be used to find exponentially many Hamilton cycles in broader classes of cubic planar graphs -- it suffices that a graph admits a nice ear decomposition and it satisfies the parity condition.

For leapfrog fullerene graphs on $n=12k$ vertices the same method can be used to find exponentially many Hamilton paths.

The real interesting open question is whether all fullerene graphs have exponentially many Hamilton cycles.

On the other hand, there are infinite series of 3-connected cubic planar graphs with faces of size at most six which only have a constant number of Hamilton cycles: Since the $\Delta$-$Y$-operation (replacement of a vertex by a triangle) conserves the number of Hamilton cycles, it suffices to consider graphs obtained from $K_4$ (the tetrahedron) by applying this operation repeatedly.

\bigskip
{\bf Acknowledgement.}
This research was partially supported by project GrR ANR-18-CE40-0032 of French National Research Agency, by project GA17–04611S of the Czech Science Foundation and by project LO1506 of the Czech Ministry of Education, Youth and Sports.

\end{document}